\documentclass[11pt, a4paper]{amsart}
\usepackage{mathbbol}
\usepackage{amsmath, amsthm, amssymb}
\usepackage{txfonts, mathrsfs}
\usepackage[usenames]{color}
\usepackage{psfrag}
\usepackage{graphicx}

\numberwithin{equation}{section}

\newtheorem{thm}{Theorem}[section]
\newcommand{\bt}{\begin{thm}}
\newcommand{\et}{\end{thm}}

\newtheorem{cor}[thm]{Corollary}
\newcommand{\bc}{\begin{cor}}
\newcommand{\ec}{\end{cor}}

\newtheorem{lem}[thm]{Lemma}
\newcommand{\bl}{\begin{lem}}
\newcommand{\el}{\end{lem}}

\newtheorem{prop}[thm]{Proposition}
\newcommand{\bp}{\begin{prop}}
\newcommand{\ep}{\end{prop}}

\newtheorem{defn}[thm]{Definition}
\newcommand{\bd}{\begin{defn}}
\newcommand{\ed}{\end{defn}}

\newtheorem{rmrk}[thm]{Remark}
\newcommand{\br}{\begin{rmrk}}
\newcommand{\er}{\end{rmrk}}

\newtheorem{quest}[thm]{Question}
\newcommand{\bq}{\begin{quest}}
\newcommand{\eq}{\end{quest}}

\newcommand{\R}{\mathbb{R}}

\newdimen\vintkern\vintkern12pt
\def\vint{-\kern-\vintkern\int}

\newcommand{\hm}{{\mathcal H}}

\newcommand{\dist}{\operatorname{dist}}

\newcommand{\trace}{\operatorname{tr}}
\newcommand{\length}{\ell}
\newcommand{\Area}{\operatorname{Area}}
\newcommand{\J}{\mathbf{J}}

\newcommand{\md}{\operatorname{md}}

\newcommand{\ap}{\operatorname{ap}}
\newcommand{\apmd}{\ap\md}

\begin{document}

\title[Harmonic discs in metric spaces]{Regularity of harmonic discs in spaces with quadratic isoperimetric inequality}


\author{Alexander Lytchak}

\address
  {Mathematisches Institut\\ Universit\"at K\"oln\\ Weyertal 86 -- 90\\ 50931 K\"oln, Germany}
\email{alytchak@math.uni-koeln.de}

\author{Stefan Wenger}

\address
  {Department of Mathematics\\ University of Fribourg\\ Chemin du Mus\'ee 23\\ 1700 Fribourg, Switzerland}
\email{stefan.wenger@unifr.ch}

\date{\today}

\thanks{The second author was partially supported by Swiss National Science Foundation Grant 153599}

\begin{abstract}
 We  study harmonic and quasi-harmonic discs in metric spaces admitting a uniformly local quadratic isoperimetric inequality for curves. The class of such metric spaces includes compact Lipschitz manifolds, metric spaces with upper or lower curvature bounds in the sense of Alexandrov, some sub-Riemannian manifolds, and many more. In this setting, we prove local H\"older continuity and continuity up to the boundary of harmonic and quasi-harmonic discs. 
 \end{abstract}

\maketitle

\maketitle
\renewcommand{\theequation}{\arabic{section}.\arabic{equation}}
\pagenumbering{arabic}

\section{Introduction and statement of main results}\label{sec:Intro}
\subsection{Background}
Questions around the existence and regularity of energy minimizing harmonic maps in various settings have been the topic of research for many years. In his pioneering work \cite{Mor48}, Morrey proved regularity of energy minimizing harmonic maps from a two-dimensional surface to a homogeneously regular Riemannian manifold. A regularity theory for higher-dimensional energy minimizing harmonic maps was developed in \cite{SU82}.
More recently, harmonic maps with values in singular metric spaces have been introduced and studied for example in \cite{GS92} and \cite{KS93}, with a particular emphasis on metric spaces of non-positive curvature in the sense of Alexandrov. In \cite{KS93} it was proved that harmonic maps from Euclidean domains to a metric space of non-positive curvature in the sense Alexandrov are locally Lip\-schitz continuous. In \cite{MZ10} it was shown that harmonic maps from the $2$-dimensional unit disc in $\R^2$ to certain compact Alexandrov spaces are locally H\"older continuous in the interior and continuous up to the boundary. Finally, the paper \cite{CL01} establishes local Lipschitz regularity for energy minimizing harmonic maps from Euclidean domains to the Heisenberg groups endowed with a Carnot-Carath\'eodory distance. The study of harmonic maps between singular metric spaces has also recently gained momentum. For regularity results for energy minimizing harmonic (real-valued) functions and, more generally, quasi-harmonic functions defined on a metric space domain we refer to \cite{KS01}, \cite{KRS03}. Regularity results for harmonic maps defined on certain singular metric spaces and with values in metric spaces of non-positive curvature can be found in \cite{Jos97}, \cite{Lin97}, \cite{Fug05}, \cite{Fug08}, \cite{MD10}, \cite{ZZ14}.

In the present paper we consider harmonic and quasi-harmonic maps from two-dimensional Euclidean domains to metric spaces. The aim of the paper is to prove H\"older regularity for such maps under much weaker conditions on the target space than was considered in the papers cited above.

\subsection{Main regularity results}
There exist several equivalent definitions of Sobolev maps defined on a Euclidean domain and with values in a complete metric space, see e.g.~\cite{KS93}, \cite{Jos94}, \cite{Haj96},  \cite{Res97}, \cite{Che99}, \cite{HKST01}, \cite{HKST15}, \cite{Res04}, \cite{Res06} and Section~\ref{sec:prelims} below. 

We briefly recall the definition proposed by Reshetnyak \cite{Res97}.
Let $X$ be a complete metric space, $\Omega\subset\R^2$ a bounded open set, and $p>1$. A measurable and essentially separably valued map $u\colon\Omega\to X$ is said to belong to $W^{1,p}(\Omega, X)$ if there exists $g\in L^p(\Omega)$ such that for every $1$-Lipschitz function $f\colon X\to \R$ the composition $f\circ u$ lies in the classical Sobolev space $W^{1,p} (\Omega)$ and its weak gradient satisfies
\begin{equation}\label{eq:Resh-upp-bd-gradient}
 |\nabla(f\circ u)|\leq g
\end{equation} almost everywhere in $\Omega$.
The Reshetnyak $p$-energy of $u\in W^{1,p}(\Omega, X)$, denoted by $E_+^p(u)$, is defined to be the infimum of the integral of $g^p$ over all $g$ satisfying \eqref{eq:Resh-upp-bd-gradient}. 

It is well-known that $u\in W^{1,p}(\Omega, X)$ in the sense above if and only if $u$ has a representative belonging to the Newton-Sobolev space $N^{1,p}(\Omega, X)$ in the sense of \cite{HKST15}. It can be shown that the Reshetnyak energy $E_+^p(u)$ equals the $p$-th power of the $L^p$-norm of the minimal weak upper gradient of such a representative of $u$ in the Newton-Sobolev space $N^{1,p}(\Omega, X)$.
 We note that a different energy was introduced by Korevaar-Schoen in \cite{KS93} which generalizes the classical  Dirichlet energy.
Finally, recall that every map $u\in W^{1,p} (\Omega,X)$ defined on a bounded Lipschitz domain $\Omega\subset\R^2$ has a well-defined trace $\trace(u)$ which is in $L^p (\partial \Omega,X)$, see \cite{KS93} and Section~\ref{sec:prelims} below.

In this paper we consider harmonic maps and, more generally, quasi-harmonic maps defined as follows.

\bd
Let $\Omega\subset\R^2$ be a bounded Lipschitz domain. A map $u\in W^{1,2}(\Omega,X)$ is called $M$-quasi-harmonic if $$E_+^2 (u|_{\Omega'})\leq M\cdot E^2_+(v)$$ for every Lipschitz domain $\Omega'\subset \Omega$ and every $v\in W^{1,2}(\Omega', X)$ with $\trace(v)= \trace(u|_{\Omega'})$.
\ed

Note that the class of quasi-harmonic maps remains unchanged if the Reshetnyak energy is replaced for example by Korevaar-Schoen's energy mentioned above or by any other definition of energy in the sense of \cite{LW-energy-area}. Harmonic maps are particular examples of $1$-quasi-harmonic maps. Here, a map $u\in W^{1,2}(\Omega, X)$ is said to be harmonic if $u$ has minimal $E_+^2$-energy among all maps in $W^{1,2} (\Omega,X)$ with the same trace as $u$.  If $X$ is a proper metric space (closed, bounded subsets of $X$ are compact) and $v\in W^{1,2}(\Omega, X)$ is such that its trace is essentially contained in a bounded ball then there exists a harmonic map with the same trace as that of $v$, see Theorem~\ref{thm:harmonic-existence}.
The class of quasi-harmonic maps has several useful properties. The restriction of an $M$-quasi-harmonic map $u\in W^{1,2} (\Omega,X)$ to a Lipschitz subdomain of $\Omega$ is again
$M$-quasi-harmonic. The class of quasi-harmonic maps is also invariant under biLipschitz changes of the metric on $X$ or on $\Omega$.

The main class of target spaces considered in the present paper are spaces admitting a uniformly local quadratic isoperimetric inequality for curves. In what follows, the open Euclidean unit disc in $\R^2$ will be denoted by $D$. We refer to Section~\ref{sec:prelims} for the precise definition of the parametrized (Hausdorff) area $\Area(u)$ of a Sobolev map $u\in W^{1,2}(\Omega, X)$ used below. 

\bd
 A complete metric space $X$ is said to admit a uniformly local quadratic isoperimetric inequality if there exist $C, l_0>0$ such that every Lipschitz curve $c\colon S^1\to X$ of length $\length(c)\leq l_0$ is the trace of a map $u\in W^{1,2} (D,X)$ with $\Area(u)\leq C\cdot\length(c)^2$.
\ed

If the constants matter we will also say that $X$ admits a $(C, l_0)$-isoperimetric inequality. Many classes of metric spaces admit a uniformly local quadratic isoperimetric inequality. This includes homogeneously regular Riemannian manifolds in the sense of Morrey \cite{Mor48}, compact Lipschitz manifolds and, in particular, compact Finsler manifolds. It also includes complete ${\rm CAT}(\kappa)$-spaces for all $\kappa\in\R$, compact Alexandrov spaces, and more generally complete metric spaces all of whose balls up to a certain radius are Lipschitz contractible in the sense of \cite{Wen07}. It furthermore includes certain sub-Riemannian manifolds such as the Heisenberg groups $\mathbb{H}^n$ of topological dimension $2n+1\geq 5$, endowed with a Carnot-Carath\'eodory metric. We refer to  \cite{LW15-Plateau} for more examples. 

We turn to the statements of our main regularity results. 

\begin{thm} \label{interior}
Let $X$ be a proper metric space admitting a $(C,l_0)$-isoperimetric inequality and $\Omega\subset \R^2$ a bounded Lipschitz domain.
If $u\in W^{1,2} (\Omega,X)$ is $M$-quasi-harmonic then $u$ has a locally $\alpha$-H\"older continuous representative $\bar{u}$ for some $\alpha$ depending only on $C$ and $M$. Moreover, if $\trace(u)$ has a continuous representative then $\bar{u}$ extends continuously to the boundary $\partial \Omega$.
\end{thm}

In particular, it follows from the theorem and \cite{Zap14} that the locally H\"older continuous representative of a quasi-harmonic map satisfies Lusin's property (N).

Our second result shows that in the same setting as above a quasi-harmonic map with Lipschitz trace is globally H\"older continuous. For simplicity, we state our theorem for the unit disc only.

\begin{thm} \label{boundary+}
Let $X$ be a proper metric space admitting a $(C,l_0)$-isoperimetric inequality. Let $u\in W^{1,2} (D,X)$ be an $M$-quasi-harmonic map whose trace has a representative which is Lipschitz continuous. Then $u$ has a representative which is globally $\alpha$-H\"older continuous for some $\alpha$ depending only on $C$, $M$.
\end{thm}

We note that H\"older continuity in Theorems~\ref{interior} and \ref{boundary+} is the best one may expect, even for harmonic maps, and cannot in general be improved to local Lipschitz continuity. This is shown by simple examples of cones over small circles, see e.g.~\cite{MZ10} or \cite{LW15-Plateau}.
Theorems~\ref{interior} and \ref{boundary+} (partly) generalize the regularity results obtained in \cite{Mor48}, \cite{KS93}, \cite{Ser94}, \cite{MZ10}, \cite{CL01}.

\subsection{An energy filling inequality}
The following theorem, which can be considered as the main theorem of the present paper, is the main new ingredient in the proof of our regularity results above. As will be shown in Section~\ref{sec:proofs-reg-results}, our regularity results above follow from the following theorem together with classical arguments going back to Morrey.

\begin{thm}\label{thm:Sobolev-fillings-energybdd}
Let $X$ be a proper metric space admitting a $(C,l_0)$-isoperimetric inequality. Then for every Lipschitz curve $c\colon S^1\to X$, parametrized proportional to arc length, with $\length(c)\leq l_0$ there exists $u\in W^{1,p}(D, X)$ with $\trace(u) = c$ and such that
\begin{equation} \label{eq:energy-filling-intro}
 \left[E_+^p(u)\right]^{\frac{1}{p}}\leq C'\cdot \length(c)
\end{equation}
for some $p>2$ and $C'$ depending only on $C$.
\end{thm}

While the geometric meaning of \eqref{eq:energy-filling-intro} might not be as clear as the meaning of the quadratic isoperimetric inequality it has several advantages due to better stability properties. For example, Theorem~\ref{thm:Sobolev-fillings-energybdd} can be used to prove stability of a uniformly local quadratic isoperimetric inequality under pointed Gromov-Hausdorff convergence of metric spaces, see \cite{LW15-asymptotic}.

Note that the map $u$ in Theorem~\ref{thm:Sobolev-fillings-energybdd} has several useful properties. Firstly, it provides an isoperimetric filling of $c$ in the sense that $\Area(u)\leq C''\cdot \length(c)^2$ for some constant $C''$ depending only on $C$. Secondly, $u$ has a representative $\bar{u}$ which is H\"older continuous on all of $D$, that is,
$$
d(\bar{u}(z), \bar{u}(z')) \leq L\cdot \nu\cdot |z-z'|^{1-\frac{2}{p}}
$$
for all $z,z'\in D$, where $\nu$ is the Lipschitz constant of $c$ and $L$ depends only on $C$, see e.g.~\cite[Proposition 3.3]{LW15-Plateau}. We furthermore note that any complete metric space $X$ in which the conclusion of Theorem~\ref{thm:Sobolev-fillings-energybdd} holds admits a $(C,l_0)$-isoperimetric inequality with $C=\frac{C'}{2\pi}$.

Theorem~\ref{thm:Sobolev-fillings-energybdd} together with a simple reparametrization argument implies the following energy filling inequality, which is prominently used in the proof of the regularity results stated above, see Section~\ref{sec:proofs-reg-results} . If $X$ is a proper metric space admitting a $(C,l_0)$-isoperimetric inequality then there exists $C'$ depending only on the isoperimetric constant $C$ such that every continuous $c\in W^{1,2}(S^1, X)$ with $\length(c)\leq l_0$ is the trace of a map $u\in W^{1,2}(D, X)$ satisfying 
\begin{equation}\label{eq:energy-filling-concl-intro}
 E_+^2(u)\leq C'\cdot E^2(c),
\end{equation}
where $E^2(c)$ denotes the energy of the curve $c$.
The filling energy inequality \eqref{eq:energy-filling-concl-intro} combined with classical arguments due to Morrey yield Theorem~\ref{interior} with H\"older exponent $\alpha=\frac{1}{2MC'}$ and Theorem~\ref{boundary+} with $\alpha=\frac{1}{\lambda MC'}$, where $\lambda$ is a universal constant. The constant $C'$ in \eqref{eq:energy-filling-concl-intro} which we obtain, and thus also the H\"older exponents $\alpha$, are far from optimal. Much more refined methods than the one used in our proofs would be needed in order to obtain optimal regularity results. In \cite{LW15-Plateau} and \cite{LW-energy-area} we obtained optimal H\"older exponents for maps which are parametrized quasi-conformally in the sense of \cite{LW15-Plateau} and minimize area rather than energy among all maps with the same trace. Note that such maps $u$ are quasi-harmonic since every $v\in W^{1,2}(\Omega')$ with $\trace(v) = \trace(u|_{\Omega'})$ satisfies
$$
 E_+^2(u|_{\Omega'}) \leq Q^2\cdot\Area(u|_{\Omega'}) \leq Q^2\cdot \Area(v) \leq Q^2 \cdot E_+^2(v),
$$ 
where $Q$ is the quasi-conformality factor of $u$. Thus, our results here generalize qualitative properties of quasi-conformal area minimizers to the class of quasi-harmonic maps.

\subsection{Outline of proof and questions} 
The proof of Thereom~\ref{thm:Sobolev-fillings-energybdd} relies on the existence results for solutions to Plateau's problem in the setting of proper metric spaces established in \cite{LW15-Plateau}.  We briefly outline the proof and let $c$ be a Lipschitz curve as in the theorem. After possibly rescaling the metric on $X$ we may assume that $c$ has length $2\pi$. Consider the enlarged space $Y=X\times\R^2$ and the uniformly biLipschitz curve $\hat{c}(z):= (c(z), z)$ in $Y$. We now solve Plateau's problem for the curve $\hat{c}$ using the results in \cite{LW15-Plateau} and obtain an area minimizing map $v\in W^{1,2}(D, Y)$ whose trace is a reparametrization of $\hat{c}$ and which is moreover $\sqrt{2}$-quasi-conformal in the sense of \cite{LW15-Plateau}. Since $Y$ admits a uniformly local quadratic isoperimetric inequality (with some constant only depending on $C$) it follows from the area minimizing property and the quasi-conformality that $v$ has $E_+^2$-energy bounded above by a constant multiple of $\length(\hat{c})^2$. Moreover, $v$ satisfies a "global" weak reverse H\"older inequality and hence $v\in W^{1,p}(\Omega, Y)$ for some $p>2$ and $E_+^p(v)^{\frac{1}{p}}\leq C''E_+^2(v)^{\frac{1}{2}}$, see Theorem~\ref{thm:global-higher-regularity}. Finally, the biLipschitz property of $\hat{c}$ can be used to find a Sobolev annulus with suitably bounded energy that reparametrizes $\trace(v)$ to $\hat{c}$, thus yielding a Sobolev map with trace $\hat{c}$. Projecting this map to $X$ yields the desired map $u$ in Theorem~\ref{thm:Sobolev-fillings-energybdd}.

As mentioned above, the constant $C'$ which we get in the theorem is far from optimal and thus yields non-optimal H\"older exponents in the regularity results. A first step towards obtaining more efficient bounds would be to understand the optimal isoperimetric constant for products of metric spaces. We formulate this as a question, compare with Lemma~\ref{lem:isop-product}.

\begin{quest}
Let $X$ be a proper geodesic metric space which is not a tree and admits a quadratic isoperimetric
inequality with constant $C$. Does the direct product $X\times \R$ admit a quadratic isoperimetric
inequality with the same constant $C$?
\end{quest}

We would also like to emphasize that properness of the underlying space is crucial in our arguments but we do not know to which extent this is actually essential for the conclusions of the theorems above.

The structure of the paper is as follows. In Section~\ref{sec:prelims} we fix notation and recall some basic definitions from the theory of Sobolev maps from a Euclidean domain to a complete metric space, including the definition of parametrized Hausdorff area for a Sobolev disc. In Section~\ref{sec:isop-implies-energy} we first establish global higher regularity of quasi-conformal area minimizers whose trace parametrizes a chord-arc curve, see Theorem~\ref{thm:global-higher-regularity}. We use this in order to prove our main result, Theorem~\ref{thm:Sobolev-fillings-energybdd}. Finally, Section~\ref{sec:proofs-reg-results} contains the proofs of our regularity results for quasi-harmonic discs.

\section{Preliminaries}\label{sec:prelims}
\subsection{Notation}
The Euclidean norm of a vector $v\in\R^n$ will be denoted by $|v|$. The unit sphere in $\R^n$ with respect to the Euclidean norm is denoted by $S^{n-1}$ and will usually be endowed with the Euclidean metric. The open unit disc in Euclidean $\R^2$ will be denoted by $D$. For $r > 0$  we denote by $\gamma _r \colon S^1\to \R^2$ the
curve $$\gamma _r (z)= r\cdot z.$$

Let $(X,d)$ be a metric space and $x\in X$. The open ball in $X$ centered at $x$ of radius $r>0$ will be denoted by $$B(x,r):= \{x'\in X: d(x,x')<r\}.$$ 
%
The length of a curve $c\colon I\to X$ from an interval $I\subset\R$ is defined by $$\length(c):= \sup\left\{\sum_{k=0}^{m-1} d(c(t_k), c(t_{k+1})):  \text{$t_0<t_1<\dots<t_m$ and $t_i\in I$}\right\}.$$ This definition readily extends to curves defined on $S^1$.
A map $u\colon \Omega\to X$, where $\Omega\subset\R^n$, is called $(L,\alpha)$-H\"older continuous if 
\begin{equation}\label{eq:Kalpha-Hoelder}
 d(u(z), u(z')) \leq L\cdot |z-z'|^\alpha
\end{equation} for all $z,z'\in\Omega$. If \eqref{eq:Kalpha-Hoelder} only holds locally around every point then $u$ is called locally $(L,\alpha)$-H\"older continuous. 

Integration on measurable subsets of $\R^2$ will always be performed with respect to the Lebesgue measure, unless otherwise stated. If $f$ is an integrable function on $\R^2$ and $B\subset\R^2$ a measurable set of strictly positive Lebesgue measure $|B|$ then $$\vint_B f(z)\,dz:= \frac{1}{|B|}\cdot\int_Bf(z)\,dz$$ denotes the averaged integral. The Hausdorff $n$-measure on a metric space will be denoted by $\hm^n$. The normalizing constant will be chosen so that $\hm^n$ coincides with the Lebesgue measure on Euclidean $\R^n$.

\subsection{Sobolev maps with values in metric spaces}\label{sec:Sobolev-maps-metric}
There exist several equivalent definitions of Sobolev maps from a Euclidean domain with values in a metric space, see e.g.~\cite{KS93}, \cite{Jos94}, \cite{Haj96},  \cite{Res97}, \cite{Che99}, \cite{HKST01}, \cite{HKST15}, \cite{Res04}, \cite{Res06}, \cite{Chi07}.
We will use the following definition due to Reshetnyak \cite{Res97}. Let $(X,d)$ be a complete metric space, $\Omega\subset\R^n$ a bounded, open set, and $p>1$. We will only need the cases $n=1$ and $n=2$ in the present paper.

\bd
A measurable and essentially separably valued map $u\colon\Omega\to X$ belongs to the Sobolev space $W^{1,p}(\Omega, X)$ if the following properties hold:
 \begin{enumerate}
  \item for every $x\in X$ the function $u_x(z):= d(x, u(z))$ belongs to the classical Sobolev space $W^{1,p}(\Omega)$.
  \item there exists $g\in L^p(\Omega)$ such that for every $x\in X$ we have $|\nabla u_x|\leq g$ almost everywhere on $\Omega$.
 \end{enumerate}
\ed

Using a local biLipschitz homeomorphism from $\R$ to $S^1$ one defines the space $W^{1,p}(S^1, X)$. Similarly, one defines the space $W^{1,p}(S^1\times(0,1), X)$.  

Let $u\in W^{1,p}(\Omega, X)$. Then for almost every $z\in \Omega$ there exists a unique semi-norm on $\R^n$, denoted by $\apmd u_z$ and called the approximate metric derivative of $u$, such that 
$$\ap \lim _{y\to z}  \frac {d(u(y),u(z))- \apmd u_z(y-z)} {|y-z|} =0,$$ see \cite{Kar07} and \cite{LW15-Plateau}. Here, $\ap\lim$ denotes the approximate limit, see \cite{EG92}. For Sobolev maps defined on an interval or on $S^1$ we will write $|\dot{c}|(s)$ or $|c'|(s)$ instead of $\apmd c_s(1)$.

 The following notion of energy was introduced by Reshetnyak, see \cite{Res97} and \cite{LW15-Plateau}.

\bd
The Reshetnyak $p$-energy $E_+^p(u)$ of a map $u\in W^{1,p}(\Omega, X)$ is defined by $$E_+^p(u)= \int_\Omega \mathcal{I}_+^p(\apmd u_z)\,dz,$$ where, for a semi-norm $s$ on $\R^n$, we have set $\mathcal{I}_+^p(s):= \max\{s(v)^p: v\in S^{n-1}\}$.
\ed

 It can be shown that the function $$g_u(z):= \mathcal{I}_+^1(\apmd u_z)$$ is the minimal weak upper gradient (of a Newtonian representative) of $u$, see \cite[Theorem 7.1.20]{HKST15}. It thus follows that $E_+^p(u)$ is the $L^p$-norm to the power $p$ of the minimal weak upper gradient of $u$. For Sobolev maps $c$ defined on an interval $(a,b)$ the Reshetnyak $p$-energy is simply given by $$\int_a^b|\dot{c}|^p(t)\,dt$$ and will be denoted by $E^p(c)$ instead of $E_+^p(c)$. 
 
 In \cite{KS93}, Korevaar-Schoen introduced a different energy which generalizes the classical Dirichlet energy. We will not use the Korevaar-Schoen energy but mention that it agrees with the Reshetnyak energy up to a non-constant bounded multiple.

We next recall the construction of the trace of a Sobolev map in \cite{KS93}. Let $\Omega\subset\R^n$ be a bounded Lipschitz domain and $u\in W^{1,p}(\Omega, X)$. For every $z\in\partial\Omega$ there exist open sets $U\subset\R^n$ and $V\subset\R^{n-1}$ with $z\in U$ and a biLipschitz homeomorphism $\varphi\colon V\times (-1,1)\to U$ such that $\varphi(V\times (-1,0)) = U\cap \Omega$ and $\varphi(V\times\{0\}) = U\cap \partial\Omega$. For almost every $v\in V$ the map $t\mapsto u\circ\varphi(v,t)$ is in $W^{1,p}((-1,0), X)$ and hence has an absolutely continuous representative, denoted by $\bar{u}(v,t)$. For $\hm^{n-1}$-almost every point $z\in U\cap\partial\Omega$ the trace of $u$ at $z$ is defined by $$\trace(u)(z):= \lim_{t\to 0^-}\bar{u}(v,t),$$ where $v\in V$ is such that $\varphi(v,0)=z$. It was shown in \cite[Lemma 1.12.1]{KS93} that the definition of $\trace(u)$ is independent of the choice of $\varphi$. Thus, covering $\partial \Omega$ by a finite number of images of biLipschitz maps, one can define $\trace(u)$ almost everywhere on $\partial\Omega$. By \cite[Theorem 1.12.2]{KS93}, the trace $\trace(u)$ belongs to the space $L^p(\partial \Omega, X)$ of measurable and essentially separably valued maps from $\partial \Omega$ to $X$ such that the composition with the distance function to any point in $X$ is in the classical space $L^p(\partial \Omega)$.
We will need the following gluing construction given in \cite[Theorem 1.12.3]{KS93}. Let $\Omega\subset\R^n$ be a bounded Lipschitz domain and suppose $\Omega$ is the disjoint union of two Lipschitz domains $\Omega_1$, $\Omega_2$ and the Lipschitz boundary $S:=\partial \Omega_1\cap \partial \Omega_2$. Let $u_i\in W^{1,p}(\Omega_i, X)$, $i=1,2$, be such that $\trace(u_1) = \trace(u_2)$ on $S$. Then the map $u$ given by $u_i$ on $\Omega_i$ for $i=1,2$ is in $W^{1,p}(\Omega, X)$ and its energy is given by $$E_+^p(u) = E_+^p(u_1) + E_+^p(u_2).$$

We have the following existence result for energy minimizing maps.

\bt\label{thm:harmonic-existence}
 Let $X$ be a proper metric space and $\Omega\subset\R^n$ a bounded Lipschitz domain. Let $p>1$ and $v\in W^{1,p}(\Omega, X)$. If there exist $x_0\in X$ and $R>0$ such that $\trace(v)(z)\in B(x_0,R)$ for almost every $z\in \partial \Omega$ then there exists $u\in W^{1,p}(\Omega, X)$ with $$E_+^p(u) = \inf\left\{E_+^p(w): w\in W^{1,p}(\Omega, X), \trace(w) = \trace(v)\right\}$$ and such that $\trace(u) = \trace(v)$.
\et

In particular, if $n=p=2$ one obtains the existence of harmonic maps with prescribed trace.

\begin{proof}
 Let $(u_j)\subset W^{1,p}(\Omega, X)$ be a minimizing sequence with respect to the Reshetnyak $p$-energy with $\trace(u_j)= \trace(v)$ for all $j$. By \cite[Lemma 3.4]{LW15-Plateau} we have $$\int_{\Omega} d^p(u_j(z), x_0)\,dz \leq L\left(R^p + E_+^p(u_j)\right)$$ for every $j$ and some $L$ depending only on $n$, $p$, $\Omega$, where $d$ denotes the metric on $X$. By Theorems 1.12.2 and 1.13 in \cite{KS93} there exists a subsequence $(u_{j_i})$ converging in $L^p(\Omega, X)$ to some $u\in W^{1,2}(\Omega, X)$ with $\trace(u) = \trace(v)$. By the weak lower semi-continuity of $E_+^p$, see \cite{Res97} or \cite[Corollary 5.7]{LW15-Plateau}, it follows that $$E_+^p(u)\leq \lim_{i\to \infty} E_+^p(u_{j_i}).$$ Thus, $u$ has minimal energy among all maps with the same trace as $v$.
\end{proof}

We now specialize to the case $n=2$. There exist several natural but different definitions of area in the literature, see \cite{LW15-Plateau}. Throughout this article, we will only work with the Hausdorff area defined as follows. Let $\Omega\subset\R^2$ be an open and bounded set in the plane. 

\bd
 The parametrized (Hausdorff) area of a map $u\in W^{1,2}(\Omega, X)$ is defined by $$\Area(u)= \int_\Omega \J(\apmd u_z)\,dz,$$ where the Jacobian $\J(s)$ of a semi-norm $s$ on $\R^2$ is the Hausdorff $2$-measure on $(\R^2, s)$ of the unit square if $s$ is a norm and zero otherwise.
 \ed
  If $u$ is injective and satisfies Lusin's property (N) then $\Area(u) = \hm^2(u(\Omega))$ by the area formula \cite{Kir94}, \cite{Kar07}.  By \cite[Lemma 7.2]{LW15-Plateau}, we have $$\J(\apmd u_z)\leq (g_u(z))^2$$ for almost every $z\in\Omega$. 

We will need the following infinitesimal notion of quasi-conformality from \cite{LW15-Plateau}. A semi-norm $s$ on $\R^2$ is called $Q$-quasi-conformal if $s(v)\leq Q\cdot s(w)$  for all $v,w\in S^1$. A map $u\in W^{1,2}(\Omega, X)$ is called $Q$-quasi-conformal if $\apmd u_z$ is $Q$-quasi-conformal for almost every $z\in\Omega$. Note that this is a much weaker notion of quasi-conformality than the one used in the field of analysis on metric spaces since $u$ is not required to be a homeomorphism. If $u$ is $Q$-quasi-conformal then 
\begin{equation}\label{eq:upper-grad-qc-jacobian}
 Q^{-2}\cdot (g_u(z))^2\leq \J(\apmd u_z)
\end{equation}
 for almost every $z\in\Omega$, see \cite[Lemma 7.2]{LW15-Plateau}.

We will need the following simple energy estimate.

\bl\label{lem:restriction-to-curves-energy}
 Let $p>1$ and let $u\in W^{1,p}(D, X)$. Then for every $s\in(0,1)$ there exists $t\in(s, 1)$ such that $u\circ\gamma_t\in W^{1,p}(S^1, X)$ with $$E^p(u\circ\gamma_t)\leq \frac{1}{1-s}\cdot E_+^p(u),$$ where $\gamma_t$ is the closed curve defined by $\gamma _t (z)= t\cdot z$ for all $z\in S^1$.
\el

\begin{proof}
 By \cite[Proposition 4.10]{LW15-Plateau} and its proof, we have $u\circ\gamma_r\in W^{1,p}(S^1, X)$ with $$|(u\circ\gamma_r)'|(v) = \apmd u_{\gamma_r(v)}(\gamma'_r(v))$$ for almost every $r\in(0,1)$ and $v\in S^1$. In particular, $$|(u\circ\gamma_r)'|^p(v) \leq r^p\cdot \mathcal{I}_+^p(\apmd u_{\gamma_r(v)})$$ for such $r$ and $v$.
 Integration in polar coordinates thus yields
 \begin{equation*}
  E_+^p(u) = \int_0^1\int_{S^1}r\cdot \mathcal{I}_+^p(\apmd u_{\gamma_r(v)})\,d\hm^1(v)\,dr \geq \int_0^1 r^{1-p} \cdot E^p(u\circ\gamma_r)\,dr,
 \end{equation*}
 from which the claim follows.
\end{proof}

\subsection{Curves and reparametrizations}\label{sec:curves-reparam}

Let $I$ denote the open unit interval in $\R$ or the unit circle $S^1$ and let $c\in W^{1,p}(I, X)$ for some $p>1$ and some complete metric space $X$. Then the continuous representative of $c$, denoted by $c$ again, is an absolutely continuous and thus rectifiable curve, which can be extended continuously to $\bar{I}$. Let $\bar{c}$ be the constant speed parametrization of $c$, see e.g.~\cite[Proposition 2.5.9]{BBI01}. We will need the following elementary observation.

\begin{lem} \label{lem:param-circle}
 There exists a homotopy $\varphi\colon \bar{I}\times [0,1]\to X$ from $c$ to $\bar{c}$ relative to endpoints which belongs to $W^{1,p}(I\times(0,1), X)$ and satisfies $\Area(\varphi)=0$ and $E_+^p(\varphi) \leq M\cdot E^p(c)$, where $M$ depends only on $p$.
\end{lem}

If $I=S^1$ then {\it relative to endpoints} should mean that $\varphi(1, t) = c(1)$ for all $t$.

\begin{proof}
 We may assume $I$ to be the unit interval, the case of the circle being almost identical. We may furthermore assume that $\lambda:=\ell(c)>0$. The normalized length function $\varrho\colon[0, 1]\to \R$ defined by  
 \begin{equation*}
 \varrho(s)= \lambda^{-1}\cdot \length(c|_{[0,s]})
\end{equation*}
satisfies $c=\bar{c}\circ\varrho$. Moreover,  $\varrho\in W^{1,p}(I)$ with $\varrho'(s) = \lambda^{-1} \cdot |\dot{c}|(s)$ for almost every $s\in[0,1]$. Define $$\psi(s,t)=(1-t)\cdot s + t\cdot \varrho(s)$$ for all $s\in\bar{I}$ and all $t\in [0,1]$. Then $\varphi:= \bar{c}\circ\psi$ is a homotopy from $c$ to $\bar{c}$ relative to endpoints. Since $\psi$ is in $W^{1,p}$ it follows that $\varphi\in W^{1,p}(I\times(0,1), X)$. Clearly, $\varphi$ has zero area and a direct calculation shows that 
$$E_+^p(\varphi) \leq M\lambda^p\cdot \int_0^1|\varrho'(s)|^p\,ds = M\cdot E^p(c)$$ for some $M$ depending only on $p$.
\end{proof}

\section{Isoperimetric inequality implies energy filling inequality}\label{sec:isop-implies-energy}

The main aim of this section is to prove Theorem~\ref{thm:Sobolev-fillings-energybdd} stated in the introduction. We will first prove a global higher integrability  result for quasi-conformal area minimizers with chord-arc boundary.

\subsection{Global higher integrability}
Let $(X,d)$ be a complete metric space and $\Gamma\subset X$ a Jordan curve. We denote by $\Lambda(\Gamma, X)$  the family of Sobolev maps $u\in W^{1,2}(D, X)$ whose trace $\trace(u)$ has a representative which weakly monotonically parametrizes $\Gamma$. This means that a representative of $\trace(u)$ is the uniform limit of homeomorphsims from $S^1$ to $\Gamma$. Recall that a Jordan curve $\Gamma\subset X$ is called $\lambda$-chord-arc curve if for any $x, y\in\Gamma$ the length of the shorter of the two segments in $\Gamma$ connecting $x$ and $y$ is bounded from above by $\lambda\cdot d(x,y)$.

\bt\label{thm:global-higher-regularity}
 Let $X$ be a complete metric space admitting a $(C,l_0$)-isoperimetric inequality, and let $\Gamma\subset X$ be a $\lambda$-chord-arc curve. If $u\in\Lambda(\Gamma, X)$ is $Q$-quasi-conformal and satisfies $$\Area(u) = \inf\{\Area(v): v\in\Lambda(\Gamma, X)\}$$ then $u\in W^{1,p}(D, X)$ for some $p>2$ depending only on $C$, $Q$, $\lambda$. 
 \et 

It follows, in particular, that $u$ has a representative which is globally $\alpha$-H\"older continuous with $\alpha=1-\frac{2}{p}$, see \cite[Proposition 3.3]{LW15-Plateau}. The proof of the theorem will furthermore show that if $\Area(u)\leq Cl_0^2$ then, after possibly precomposing $u$ with a Moebius transformation, we have 
 \begin{equation}
  E_+^p(u)\leq L\cdot \left[E_+^2(u)\right]^{\frac{p}{2}}
 \end{equation}
  for some constant $L$ depending only on $C$, $Q$, $\lambda$.

The proof of Theorem~\ref{thm:global-higher-regularity} can be obtained by combining the arguments in the proofs of Theorems 8.2 and 9.3 in \cite{LW15-Plateau}. For the sake of completeness we sketch the argument.

\begin{proof}
 By \cite[Theorem 1.4]{LW15-Plateau} we may assume that $u$ is continuous on all of $\overline{D}$.
 Fix three points $p_1, p_2, p_3\in S^1$ at equal distance from each other and let $q_1,q_2,q_3\in\Gamma$ be three points such that the three segments into which they divide $\Gamma$ have equal length. After possibly precomposing $u$ with a Moebius transformation we may assume that $u$ satisfies the $3$-point condition $u(p_i) = q_i$ for $i=1,2,3$.
 
Throughout the proof, we will denote by $C_i$, $i=1,2,\dots$, constants that only depend on $C$, $Q$, and $\lambda$. Let $0<r_0\leq \frac{1}{4}$ be such that $\Area(u|_{D\cap B(z,2r_0)})\leq Cl_0^2$ for every $z\in \overline{D}$. Define a non-negative function $f\in L^2(\R^2)$ by $f(z)= \mathcal{I}_+^1(\apmd u_z)$ if $z\in D$ and $f(z)= 0$ otherwise. We first show that for every square $W\subset\R^2$ of side length at most $2r_0$ the weak reverse H\"older inequality 
\begin{equation}\label{eq:reverse-Hoelder}
  \left(\,\vint_W f^2(z)\,dz\right)^{\frac{1}{2}} \leq C_1\cdot \vint_{2W}f(z)\,dz
\end{equation}
holds, where $2W$ denotes the square with the same center as $W$ but twice the side length. Let $W\subset\R^2$ be a square of side length $2r\leq 2r_0$ centered at some point $z\in\R^2$. We may assume that $W\cap D\not=\emptyset$. The proof of \cite[Theorem 9.3]{LW15-Plateau} shows that $$\Area(u|_{D\cap B(z,s)})\leq C(1+2\lambda)^2\cdot \length(u|_{D\cap \partial B(z,s)})^2$$ for almost every $s\in(\sqrt{2}r, 2r)$. Since $$\length(u|_{D\cap \partial B(z,s)})\leq \int_{D\cap \partial B(z,s)} f(z)\,d\hm^1(z)$$ for almost every $s$, see e.g.~the proof of \cite[Proposition 8.4]{LW15-Plateau}, we obtain together with \eqref{eq:upper-grad-qc-jacobian} that
\begin{equation*}
 \begin{split}
   \left(\int_Wf^2(z)\,dz\right)^{\frac{1}{2}} &\leq Q\cdot \Area(u|_{D\cap W})^{\frac{1}{2}}\\
  &\leq Q\sqrt{C} (1+2\lambda)\cdot\vint_{\sqrt{2}r}^{2r} \int_{D\cap \partial B(z,s)} f(z)\,d\hm^1(z)\,ds\\
  &\leq \frac{Q\sqrt{C}(1+2\lambda)}{(2-\sqrt{2}) r}\cdot \int_{2W}f(z)\,dz.
 \end{split}
\end{equation*}
This implies \eqref{eq:reverse-Hoelder}. By the generalized Gehring lemma, see e.g.~\cite[p.~409]{Str80} or \cite[Theorem 1.5]{Kin94}, there thus exists $p>2$ depending only on $C_1$, and thus only on $C$, $Q$, $\lambda$, such that 
\begin{equation}\label{eq:implication-Gehring}
 \left(\,\vint_W f^p(z)\,dz\right)^{\frac{1}{p}} \leq C_2\cdot \left(\,\vint_{2W}f^2(z)\,dz\right)^{\frac{1}{2}}
\end{equation}
 for every square $W\subset\R^2$ of side length at most $2r_0$. Covering $D$ by almost disjoint squares of side length $2r_0$ and using \eqref{eq:implication-Gehring} one obtains $$\int_D f^p(z)\,dz \leq C_3 r_0^{2-p}\cdot\left(\int_Df^2(z)\,dz\right)^{\frac{p}{2}}$$ and hence $f\in L^p(\R^2)$. It thus follows that $u\in W^{1,p}(D, X)$ with $$E_+^p(u)\leq C_3 r_0^{2-p}\cdot\left[ E_+^2(u)\right]^{\frac{p}{2}}.$$ This completes the proof.
\end{proof}

\subsection{Proof of Theorem~\ref{thm:Sobolev-fillings-energybdd}}\label{sec:Proof-Sob-filling-energybdd}

The following simple construction, which produces a new Sobolev disc from an old one by attaching a Sobolev annulus, will be employed several times in the proofs below.
Let $v\in W^{1,p}(D, X)$ be a Sobolev disc for some $p>1$. Let $\Omega\subset\R^2$ be an annulus of the form $\Omega=\{1<|z|<t\}$ for some $t\in(1, 2]$, and let $w\in W^{1,p}(\Omega, X)$. If $\trace(v) = \trace(w)|_{S^1}$ then the gluing construction mentioned in Section~\ref{sec:Sobolev-maps-metric} lets us glue $v$ and $w$ along $S^1$ in order to obtain a new Sobolev map defined on the disc of radius $t$. Thus, after rescaling the domain, we obtain a map $v'\in W^{1,p}(D, X)$ with $\trace(v') = \trace(w)\circ\gamma_t$ and such that $$\Area(v') = \Area(v) + \Area(w)$$ and 
\begin{equation}\label{eq:energy-gluing-annulus}
E_+^p(v')\leq L\cdot E_+^p(v) + L\cdot E_+^p(w)
\end{equation}
 for some $L$ depending only on $p$. The map $v'$ is said to be obtained by attaching the Sobolev annulus $w$ to the Sobolev disc $v$.
By possibly precomposing with a biLipschitz homeomorphism one can also attach a Sobolev annulus defined on any set biLip\-schitz homeomorphic to $\Omega$ such as for example $S^1\times (0,1)$. In this case the constant $L$ in \eqref{eq:energy-gluing-annulus} also depends on the biLip\-schitz constant.

\bl\label{lem:isop-product}
 If a complete metric space $X$ admits a $(C,l_0)$-isoperimetric inequality then the product space $Y=X \times \R^2$ admits a $(C',l_0)$-isoperimetric inequality with $C'=C+M$ for some universal constant $M$.
\el

The proof will in fact show that every Lipschitz curve $c$ in $Y$ with the property that its projection to $X$ has length bounded by $l_0$ is the trace of some $u\in W^{1,2}(D, Y)$ with $\Area(u) \leq C'\cdot \length(c)^2$.

\begin{proof}
 Let $c=(c_1,c_2)\colon S^1\to Y$ be a Lipschitz curve such that the component $c_1$ has length at most $l_0$. In view of Lemma~\ref{lem:param-circle} and the gluing construction above, we may assume that $c$ is parametrized proportional to arc length.
Define a curve $\hat{c}$ in $Y$ by $\hat{c} (z)= (c_1(z),c_2(1))$. Since $\hat{c}$ has length $\length(\hat{c})=\length(c_1)\leq l_0$ and since its image is contained in a subspace of $Y$ isometric to $X$ there exists a Sobolev disc $v\in W^{1,2}(D, Y)$ whose trace equals $\hat{c}$ and whose area satisfies $$\Area(v)\leq C\cdot\length(c_1)^2\leq C\cdot\length(c)^2.$$
We next define a homotopy $w\colon S^1\times [0,1]\to Y$ from $c$ to $\hat{c}$ by $$w(z,t)=(c_1(z), t c_2(1) + (1-t)c_2(z)).$$ Since $c$ is parametrized proportional to arc length it follows that $w$ is $M_1\cdot \length(c)$-Lipschitz and thus has $\Area(w)\leq M_2\cdot \length(c)^2$ for some universal constants $M_1$ and $M_2$. 
We attach $w$ to $v$ along the curve $\hat{c}$ in order to obtain a Sobolev disc $u\in W^{1,2}(D, Y)$ whose trace is $c$ and whose area is $$\Area(u) = \Area(v) + \Area(w)\leq (C+M_2)\cdot\length(c)^2.$$ This completes the proof.
\end{proof}

We finally turn to the proof of Theorem~\ref{thm:Sobolev-fillings-energybdd}.
Let $c\colon S^1\to X$ be a Lipschitz curve of length $\length(c)\leq l_0$ which is parametrized proportional to arc length. Since the conclusions of the theorem are invariant under rescaling of the metric we may assume that $\length(c) = 2\pi$ and $2\pi \leq l_0$.

Consider the product space $Y=X\times \R^2$ and let $\hat{c}\colon S^1\to Y$ be the curve given by $\hat{c}(z):= (c(z), z)$ for all $z\in S^1$. Since the natural projection from $Y$ to $X$ is $1$-Lipschitz it suffices to show that there exists a map $u\in W^{1,p}(D, Y)$ with $\trace(u) = \hat{c}$ and such that $E_+^p(u)\leq C'$ for some $p>2$ and $C'$ only depending on the isoperimetric constant $C$. 
We will construct such a map $u$ in several steps. We will first exhibit a Sobolev disc $v_1$ whose trace is only a weakly monotone reparametrization of $\hat{c}$. We will then attach suitable Sobolev annuli to $v_1$ in order to obtain a map whose trace actually equals $\hat{c}$.

In what follows, $C_1$, $C_2$, \dots will denote constants depending only on $C$. We construct the Sobolev disc $v_1$ as follows. First note that $\hat{c}$ is $\lambda$-biLipschitz for some universal $\lambda$ and hence its image $\Gamma$ is a chord-arc curve in $Y$ with universal chord-arc parameter. Lemma \ref{lem:isop-product} and the subsequent remark imply that $\Lambda (\Gamma ,Y)$ is not empty and that $$m:=\inf\{\Area(v): v\in\Lambda(\Gamma, Y)\}\leq C_1\cdot \length(\hat{c})^2\leq C_2.$$ By \cite[Theorem 1.1]{LW15-Plateau} there thus exists $v_1\in \Lambda(\Gamma, Y)$ which is $\sqrt{2}$-quasi-conformal and satisfies $\Area(v_1) = m\leq C_2$. It thus follows from Theorem~\ref{thm:global-higher-regularity} that there exist $p>2$ and $L>0$ depending only on $C$ such that $v_1\in W^{1,p}(D, Y)$ and such that, after possibly precomposing $v_1$ with a M\"obius transformation, we have $$E_+^p(v_1)\leq L\cdot \left[E_+^2(v_1)\right]^{\frac{p}{2}}\leq L\cdot\left[2\Area(v_1)\right]^{\frac{p}{2}} \leq C_3.$$ In particular, $v_1$ (has a representative which) is $(C_4, \alpha)$-H\"older on all of $\overline{D}$ with $\alpha=1-\frac{2}{p}$, see e.g.~\cite[Proposition 3.3]{LW15-Plateau}. This completes the construction of the Sobolev disc $v_1$ whose trace is a weakly monotone parametrization of $\hat{c}$.

We will next attach a Sobolev annulus $w_1$ to $v_1$ in order to obtain a Sobolev disc $v_2$ whose trace is a $p$-Sobolev curve with image in $\hat{c}$. This annulus will be obtained by restricting $v_1$ to a suitably chosen, small annulus $\{t<|z|<1\}$ and then mapping it via a Lipschitz projection onto the curve $\hat{c}$. For this, observe first that $\Gamma$ is an absolute Lipschitz neighborhood retract since $\hat{c}$ is $\lambda$-biLipschitz. In particular, there exists an $M$-Lipschitz retraction $\eta\colon N_\varepsilon(\Gamma) \to\Gamma$ of the closed $\varepsilon$-neighborhood $N_\varepsilon(\Gamma)$ of $\Gamma$ in $Y$ to $\Gamma$ for some universal constants $\varepsilon, M>0$. Since $v_1$ is $(C_4, \alpha)$-H\"older there exists $s\in(0,1)$ depending only on $C$ such that $v_1(z)\in N_\varepsilon(\Gamma)$ for every $z\in D$ with $|z|\geq s$. By Lemma~\ref{lem:restriction-to-curves-energy} there exists $t\in(s,1)$ such that $v_1\circ\gamma_t\in W^{1,p}(S^1, Y)$ and $$E^p(v_1\circ\gamma_t) \leq \frac{1}{1-s}\cdot E_+^p(v_1) \leq C_5.$$ Let $w_1$ be the Sobolev annulus  given by composing the restriction of $v_1$ to the annulus $\{t<|z|<1\}$ with the map $\eta$. Attaching $w_1$ to $v_1$ along the outer boundary $S^1$ of the annulus we obtain a map $v_2\in W^{1,p}(D, Y)$ whose energy $E_+^p(v_2)$ is bounded by some constant $C_6$ and whose trace equals $c_0:=\eta\circ v_1\circ\gamma_t$. Notice that $E^p(c_0) \leq M^p C_5\leq C_7$.

We finally attach a Sobolev annulus $\varphi$ to $v_2$ in order to obtain a map $v_3$ with trace $\hat{c}$ as follows. We first observe that $c_0$ is of the form $c_0= \hat{c} \circ \varrho$ for some map $\varrho\colon S^1\to S^1$ which is homotopic to the identity. We can now define a `linear' homotopy $\varphi$ from $c_0$ to $\hat{c}$ with $$E_+^p(\varphi)\leq C_8\cdot E^p(c_0)\leq C_9$$ as in the proof of Lemma~\ref{lem:param-circle}. Attaching the annulus $\varphi$ to the map $v_2$ along the curve $c_0$ we finally obtain a map $v_3\in W^{1,p}(D, Y)$ whose trace equals $\hat{c}$ and whose $p$-energy is bounded by some constant $C_{10}$. This completes the proof of Theorem~\ref{thm:Sobolev-fillings-energybdd}.

The last part of the proof of Theorem~\ref{thm:Sobolev-fillings-energybdd} shows the following:

\bc
 Let $Y$ be a complete metric space, $\Gamma\subset Y$ a chord-arc curve, and $u\in \Lambda(\Gamma, Y)\cap C^0(\overline{D}, Y)$. Then there exists $u_0\in\Lambda(\Gamma, Y)$ such that $\trace(u_0)$ parametrizes $\Gamma$ proportionally to arc length and $\Area(u_0) = \Area(u)$.
\ec

\begin{proof}
 The Sobolev annuli $w_1$ and $\varphi$ constructed in the proof have zero area. Thus, attaching the two annuli to $u$ yields a map $u_0$ with the desired properties. 
\end{proof}

The proof of Theorem~\ref{thm:Sobolev-fillings-energybdd} together with Lemma~\ref{lem:param-circle} moreover yields the following strengthening of the theorem:

\bt\label{thm:Sobolev-filling-energy-arb-Lip}
Let $X$ be a proper metric space admitting a $(C,l_0)$-isoperimetric inequality. Then for every $\varepsilon>0$ there exist $p>2$ and $C'$ only depending on $C$, $\varepsilon$ with the following property. Every continuous curve $c\in W^{1,p}(S^1, X)$ with $\length(c)\leq l_0$ is the trace of some $u\in W^{1,p}(D, X)$ with $\Area(u)\leq (C+\varepsilon)\cdot\length(c)^2$ and $E_+^p(u)\leq C'\cdot E^p (c)$.
\et

\begin{proof}
By Lemma~\ref{lem:param-circle}, one may assume that $c$ is parametrized proportional to arc length. After rescaling the metric, we may furthermore assume that $c$ has length $\length(c)=2\pi$. In order to obtain the better area bound one considers, for $\varepsilon>0$ fixed, the curve $\hat{c}(z):= (c(z), \varepsilon\cdot z)$ instead of the curve $\hat{c}$ defined in the proof of Theorem~\ref{thm:Sobolev-fillings-energybdd}. The image $\Gamma$ of $\hat{c}$ is a chord-arc curve in  $Y$ with chord-arc parameter depending on $\varepsilon$. The Sobolev disc constructed in the proof of Lemma~\ref{lem:isop-product} yields a filling of $\hat{c}$ whose area is bounded by  $C\cdot\length(c)^2 + M\varepsilon$ for some universal constant $M$. Thus, the quasi-conformal area minimizer $v_1\in \Lambda(\Gamma, Y)$ in the proof of Theorem~\ref{thm:Sobolev-fillings-energybdd} has the same area bound. The rest of the proof is identical to that of Theorem~\ref{thm:Sobolev-fillings-energybdd}. Since the Sobolev annuli $w_1$ and $\varphi$ in the proof have zero area the filling of $c$ has area at most $C\cdot\length(c)^2 + M\varepsilon$.
\end{proof}

\section{Regularity of quasi-harmonic maps}\label{sec:proofs-reg-results}

In this section we prove the main regularity results of the present paper stated in Theorems~\ref{interior} and \ref{boundary+}. As already mentioned in the introduction, these results follow from Theorem~\ref{thm:Sobolev-fillings-energybdd} together with classical arguments going back to Morrey, which we will briefly explain below. The main implication of Theorem~\ref{thm:Sobolev-fillings-energybdd} which we will need is the following energy filling inequality.

\subsection{Energy filling inequality}\label{subsec:energy-fill}
Let $X$ be a proper metric space admitting a $(C, l_0)$-isoperimetric inequality. Then Theorem~\ref{thm:Sobolev-fillings-energybdd}, Lemma~\ref{lem:param-circle}, and the gluing construction at the beginning of Section~\ref{sec:Proof-Sob-filling-energybdd} imply that $X$ admits a uniformly local energy filling inequality in the following sense.
There exists $C'>\frac{1}{2}$ depending only on $C$ such that every continuous $c\in W^{1,2}(S^1, X)$ with $\length(c)\leq l_0$ is the trace of a map $u\in W^{1,2}(D, X)$ satisfying 
\begin{equation}\label{eq:weak-filing-energy-ineq}
 E_+^2(u)\leq C'\cdot E^2(c).
\end{equation}  
We will say that $X$ admits a $(C', l_0)$-energy filling inequality. 
Note that conversely, if $X$ is a complete metric space admitting a $(C',l_0)$-energy filling inequality in the sense above then $X$ admits a $(C,l_0)$-isoperimetric inequality with $C=\frac{C'}{2\pi}$. 
 
\subsection{Regularity in spaces with energy filling inequality} 
Theorems~\ref{interior} and \ref{boundary+} are direct consequences of the statement in Subsection~\ref{subsec:energy-fill} and the following regularity results.

\bp\label{prop:int-hoelder-energy-filling}
 Let $X$ be a complete metric space admitting a $(C', l_0)$-energy filling inequality, and let $\Omega\subset\R^2$ be a bounded Lipschitz domain. If $u\in W^{1,2}(\Omega, X)$ is $M$-quasi-harmonic then $u$ has a continuous representative which is locally $\alpha$-H\"older continuous on $\Omega$ with $\alpha = \frac{1}{2MC'}$. 
\ep

The proof of the proposition relies on the following version for metric spaces of Morrey's growth theorem. Let $r_0>0$ be fixed. For $z\in \Omega$ set 
\begin{equation}\label{eq:def-varrho-dist}
 \varrho(z):= \min\{r_0, \dist(z, \partial \Omega)\}.
\end{equation}

Then we have:

\bl\label{lem:Morrey-growth}
 Let $X$ be a complete metric space, $K\geq 0$, and $\alpha\in(0,1)$. Suppose $u\in W^{1,2}(\Omega, X)$ satisfies 
 \begin{equation}\label{eq:growth-estimate-Morrey}
  E_+^2(u|_{B(z,r)}) \leq K\cdot \varrho(z)^{-2\alpha}\cdot r^{2\alpha}
\end{equation}
 for all $z\in\Omega$ and $r\in(0, \varrho(z))$. Then $u$ is locally $\alpha$-H\"older continuous on $\Omega$. 
\el

More precisely, the proof will show that $u$ has a continuous representative which is $(K', \alpha)$-H\"older continuous on $B(z, \varrho(z)/2)$ for every $z\in\Omega$, where $K'= L\sqrt{K}\cdot \varrho(z)^{-\alpha}$ for some $L$ only depending on $\alpha$.

\begin{proof}
Let $x\in X$. Then the composition $u_x:= d_x\circ u$ of $u$ with the distance function $d_x$ to $x$ in $X$ is in $W^{1,2}(\Omega)$ with $$|\nabla u_x(w)|\leq \mathcal{I}_+^1(\apmd u_w)$$ for almost every $w\in\Omega$. Thus, H\"older's inequality and \eqref{eq:growth-estimate-Morrey} yield
$$
 \int_{B(z,r)} |\nabla u_x(w)|\,dw \leq \sqrt{\pi}r\cdot\left(E_+^2(u|_{B(z,r)})\right)^{\frac{1}{2}}\leq \sqrt{K\pi}\cdot \varrho(z)^{-\alpha}\cdot r^{1+\alpha}
$$
for every $z\in\Omega$ and every $r\in(0,\varrho(z))$.
This together with Morrey's theorem, see e.g.~\cite[Theorem 7.19]{GT01}, implies that $u_x$ has a continuous representative which is $(K', \alpha)$-H\"older continuous on $B(z, \varrho(z)/2)$ for each $z\in\Omega$, where $K'= L\sqrt{K}\cdot \varrho(z)^{-\alpha}$ for some $L$ only depending on $\alpha$. 

Using the above for a dense sequence of $x$ in the essential image of $u$ one easily obtains that also $u$ has a continuous representative which is $(K', \alpha)$-H\"older continuous on $B(z, \varrho(z)/2)$ for each $z\in\Omega$.
\end{proof}

  \begin{proof}[Proof of Proposition~\ref{prop:int-hoelder-energy-filling}]
   Let $r_0>0$ be such that $E_+^2(u|_{\Omega\cap B(z,r_0)})\leq A$ for every $z\in\Omega$, where $A= (2\pi)^{-1} MC'l_0^2$. Let $\varrho$ be defined as in \eqref{eq:def-varrho-dist} and set 
   $\alpha=\frac{1}{2MC'}$. We will show that
  \begin{equation}\label{eq:growth-proof-prop-int-hoeld}
   E_+^2(u|_{B(z,r)})\leq  E_+^2(u)\cdot\varrho(z)^{-2\alpha}\cdot r^{2\alpha}
  \end{equation}
  for all $z\in\Omega$ and $r\in(0,\varrho(z))$. The proposition will follow from this and Lemma~\ref{lem:Morrey-growth}. 
  
  Let $z\in\Omega$ and set $r_1:= \varrho(z)$. For $r>0$ let $\gamma_{z,r}\colon S^1\to \R^2$ be the curve $\gamma_{z,r}(v):= z+rv$.
 As in the proof of Lemma~\ref{lem:restriction-to-curves-energy} we have that $u\circ\gamma_{z,r}\in W^{1,2}(S^1, X)$ and 
  \begin{equation}\label{eq:Hoelder-reg-speed-minupgrad}
   |(u\circ\gamma_{z,r})'|^2(v) \leq r^2\cdot \mathcal{I}_+^2(\apmd u_{\gamma_{z,r}(v)})
  \end{equation}
  for almost every $r\in(0,r_1)$ and $v\in S^1$. We claim that for each such $r$ we have
   \begin{equation}\label{eq:energy-ineq-proof-lochoeld}
    E_+^2(u|_{B(z,r)})\leq MC'\cdot E^2(u\circ\gamma_{z,r}).
   \end{equation}
Indeed, if the length of the absolutely continuous representative of $u\circ\gamma_{z,r}$ is smaller than $l_0$ then \eqref{eq:energy-ineq-proof-lochoeld} follows from the quasi-harmonicity of $u$ and \eqref{eq:weak-filing-energy-ineq}. If this length is larger than $l_0$ then \eqref{eq:energy-ineq-proof-lochoeld} is a consequence of H\"older's inequality and the bound $E_+^2(u|_{B(z,r_1)})\leq A$. This proves \eqref{eq:energy-ineq-proof-lochoeld}. 

Now, by \eqref{eq:Hoelder-reg-speed-minupgrad} and the change of variable formula, we have $$E^2(u\circ\gamma_{z,r}) \leq r\cdot \int_{\partial B(z,r)} \mathcal{I}_+^2(\apmd u_w)\,d\hm^1(w) = r\cdot\frac{d}{dr} E_+^2(u|_{B(z,r)})$$ for almost every $r\in(0,r_1)$ and thus
    \begin{equation*}
    E_+^2(u|_{B(z,r)})\leq MC'\cdot r \cdot\frac{d}{dr} E_+^2(u|_{B(z,r)}).
   \end{equation*}
 From this, \eqref{eq:growth-proof-prop-int-hoeld} follows upon integration. 
  \end{proof}

\br\label{rem:int-hoelder-remark}
 {\rm The proof of Proposition~\ref{prop:int-hoelder-energy-filling} does not use the full strength of the definition of quasi-harmonicity. Indeed, one only needs the quasi-harmonicity condition on balls which are compactly contained in $\Omega$. The proof in particular shows that if $E_+^2(u)\leq A$, where $A$ is as in the proof, then for every $z\in\Omega$ the continuous representative $\bar{u}$ of $u$ is $(K', \alpha)$-H\"older continuous on $B(z,\varrho/2)$, where $\varrho$ is the distance of $z$ to $\partial \Omega$ and $K' = L\cdot E_+^2(u)^{\frac{1}{2}}\cdot\varrho^{-\alpha}$ for some $L\geq 1$ only depending on $\alpha$.}
\er

The second part of Theorem~\ref{interior} follows from Subsection~\ref{subsec:energy-fill} and the following proposition.

\bp\label{prop:bdry-cont}
Let $X$ be a complete metric space admitting a uniformly local energy filling inequality, and let $\Omega\subset\R^2$ be a bounded Lipschitz domain. If $u\in W^{1,2}(\Omega, X)\cap C^0(\Omega, X)$ is quasi-harmonic and $\trace(u)$ has a continuous representative then $u$ extends continuously to $\partial \Omega$.
\ep

We sketch the proof of the proposition, which is the same as that of \cite[Theorem 9.1]{LW15-Plateau} and relies on the following lemma appearing in \cite[Lemma 9.2]{LW15-Plateau}.

\bl\label{lem:energy-growth-bdry-cont}
Let $(X,d)$ be a complete metric space and $v\in W^{1,2}(D, X)\cap C^0(D,X)$. If $v$ satisfies 
$d(\trace(v)(z), v(0))\geq \varepsilon>0$ for almost every $z\in S^1$ and $d(v(z), v(0))<\frac{\varepsilon}{2}$ for all $z\in D$ with $|z|<\frac{1}{2}$
then $E_+^2(v|_{v^{-1}(B(v(0), \varepsilon))})\geq F\varepsilon^2$ for some universal constant $F\in(0,1)$.
\el

\begin{proof}[Proof of Proposition~\ref{prop:bdry-cont}]
Since quasi-harmonicity is preserved under passing to Lipschitz subdomains and under precomposing with a biLipschitz map we may assume that $\Omega=D$. Let $C', l_0, M$ be such that $X$ admits a $(C',l_0)$-energy filling inequality and $u$ is $M$-quasi-harmonic. Set $\alpha = \frac{1}{2MC'}$ and let $A$ and $L\geq 1$ be as in Remark~\ref{rem:int-hoelder-remark}. Finally, let $F\in(0,1)$ be the constant from Lemma~\ref{lem:energy-growth-bdry-cont}.

Let $\bar{u}\colon\overline{D}\to X$ be the map which coincides with $u$ on $D$ and with the continuous representative of $\trace(u)$ on $S^1$. We must show that $\bar{u}$ is continuous. For this it is enough to prove that $\bar{u}$ is continuous at each $z\in S^1$ since $\bar{u}$ is continuous on $D$. Fix $z\in S^1$ and let $0<\varepsilon<\frac{F}{2L}$. By the continuity of the restriction $\bar{u}|_{S^1}$ and the Courant-Lebesgue lemma there exists $r>0$ such that the restriction of $\bar{u}$ to the boundary of $\Omega':= D\cap B(z,r)$ is continuous and has image contained in $B(\bar{u}(z), \varepsilon)$. See the proof of \cite[Theorem 9.1]{LW15-Plateau} for details. Moreover, this $r$ can be chosen arbitrarily small and we may therefore assume that $E_+^2(u|_{\Omega'})<\min\{A, \varepsilon^4\}$. We will show that $d(\bar{u}(z), u(x))<2\varepsilon$ for every $x\in\Omega'$, thus establishing continuity of $\bar{u}$ at $z$. Fix $x\in\Omega'$ and let $\varphi\colon D\to \Omega'$ be a conformal diffeomorphism which maps $0$ to $x$. We want to apply Lemma~\ref{lem:energy-growth-bdry-cont} to the map $v:= u\circ \varphi$. Observe first that $v\in W^{1,2}(D,X)$ with $E_+^2(v)= E_+^2(u|_{\Omega'})$ and $\trace(v)= \bar{u}|_{\partial \Omega'}\circ \varphi|_{S^1}$. Moreover, since $u$ is $M$-quasi-harmonic and $\varphi$ is a conformal diffeomorphism it follows that $v$ satisfies the $M$-quasi-harmonicity condition for every Lipschitz domain compactly contained in $D$. Thus, $v$ is $(K', \alpha)$-H\"older continuous on $B(0, \frac{1}{2})$ with $K' \leq L\cdot \varepsilon^2$, by Remark~\ref{rem:int-hoelder-remark}. In particular, $$d(v(z'), v(0))\leq L\cdot\varepsilon^2\cdot|z'|^\alpha<\frac{\varepsilon}{2}$$ for every $z'$ with $|z'|<\frac{1}{2}$. Moreover, the image of $\trace(v)$ is the same as that of $\bar{u}|_{\partial \Omega'}$ and hence contained in $B(\bar{u}(z), \varepsilon)$. Thus, we must have $d(v(0), \bar{u}(z))\leq 2\varepsilon$ since otherwise $d(\trace(v)(z), v(0))> \varepsilon$ for every $z\in S^1$ and Lemma~\ref{lem:energy-growth-bdry-cont} would then give $$F\varepsilon^2\leq E_+^2(v)< \varepsilon^4 < F\varepsilon^2,$$ a contradiction. This shows that $d(u(x), \bar{u}(z))\leq 2\varepsilon$ and thus the proof is complete.
\end{proof}

We finally prove Theorem~\ref{boundary+}. This is a consequence of the statement in Subsection~\ref{subsec:energy-fill} together with Proposition~\ref{prop:bdry-cont} and the following:

\bp\label{prop:boundary-lip-global-hoelder}
  Let $X$ be a complete metric space admitting a $(C', l_0)$-energy filling inequality. Let $u\in W^{1,2} (D,X)\cap C^0(\overline{D}, X)$ be an $M$-quasi-harmonic map. If $u|_{S^1}$ is Lipschitz then $u$ is globally $\alpha$-H\"older continuous with $\alpha=\frac{1}{\lambda C'M}$, where $\lambda$ is a universal constant.
\ep

We first recall the following global analog of Morrey's theorem, which can be obtained from the corresponding classical result as in Lemma~\ref{lem:Morrey-growth}. Let $r_0>0$ and suppose $u\in W^{1,2}(D,X)$ satisfies $$E_+^2(u|_{D\cap B(z,r)}) \leq K \cdot r^{2\alpha}$$ for every $z\in D$ and every $r\in(0,r_0)$. Then $u$ has a globally $\alpha$-H\"older continuous representative. We will furthermore need the following simple observation.

\bl\label{lem:diff-ineq}
 Let $\beta\in(0,2)$ and $A\geq 0$. Let $f\colon[0,r_0]\to [0,\infty)$ be a continuous and increasing function such that $$f(r) \leq \frac{1}{\beta}\cdot r\cdot f'(r) + Ar^2$$ for almost every $r\in[0,r_0]$. Then $$f(r)\leq \left(f(r_0)\cdot r_0^{-\beta} + \bar{A}\cdot r_0^{2-\beta}\right)\cdot r^\beta$$ for every $r\in[0,r_0]$, where $\bar{A}= \frac{A\beta}{2-\beta}$.
\el

\begin{proof}
The continuous and increasing function given by $F(r):= f(r) + \bar{A} r^2$ satisfies the differential inequality 
 \begin{equation*}
   F(r)\leq \frac{1}{\beta}\cdot r\cdot f'(r) + Ar^2 + \bar{A} r^2 =  \frac{1}{\beta}\cdot r\cdot F'(r)
 \end{equation*} 
 for almost every $r\in[0,r_0]$. Upon integration, we obtain $$F(r) \leq F(r_0)\cdot r_0^{-\beta}\cdot r^\beta$$ for every $r\in[0,r_0]$, from which the result follows.
\end{proof}

\begin{proof}[Proof of Proposition~\ref{prop:boundary-lip-global-hoelder}]
 Let $0<r_0<1$ be such that for every $z\in\overline{D}$ we have $$E_+^2(u|_{D\cap B(z,r_0)})\leq (2\pi)^{-1}\cdot MC'l_0^2.$$ Fix $z\in\overline{D}$. Define a continuous and increasing function by $f(r):= E_+^2(u|_{D\cap B(z,r)})$. In view of Lemma~\ref{lem:diff-ineq} and the global analog of Morrey's theorem mentioned above it suffices to show that 
 \begin{equation}\label{eq:diff-ineq-energy}
 f(r)\leq \lambda C'M\cdot r\cdot f'(r) + Ar^2
 \end{equation} 
 for almost every $r\in(0,r_0)$, where $\lambda\geq 1$ is a universal constant and where $A\geq0$ depends on $C'$, $M$, and on the Lipschitz constant of $u|_{S^1}$.
 
By the proof of Proposition~\ref{prop:int-hoelder-energy-filling} we have \eqref{eq:diff-ineq-energy} with $\lambda=1$ and $A=0$ for almost every $r\in(0,r_0)$ such that $B(z,r)\subset D$. Let now $r\in(0,r_0)$ be such that $B(z,r)\not\subset \overline{D}$. We denote by $\gamma_r$ the curve $\gamma_r(v):= z+rv$ for $v\in S^1$. Set $\Gamma_1:= \gamma_r^{-1}(D)$ and let $\varphi$ be a biLipschitz homeomorphism from  $\overline{B}(z,r)$ to $\overline{D}\cap\overline{B}(z,r)$ with a universal biLipschitz constant $\eta$ and such that $\varphi$ is the identity on $\Gamma_1$. Such a map $\varphi$ exists by \cite{Tuk80}. 
Firstly, we claim
\begin{equation}\label{eq:fr-qh-ef}
 f(r) \leq \eta^4 C'M\cdot E^2(u\circ\varphi\circ\gamma_r).
\end{equation}
 Indeed, if $\length(u\circ\varphi\circ\gamma_r)\leq l_0$ then this is a consequence of the energy filling inequality and the quasi-harmonicity of $u$. If $\length(u\circ\varphi\circ\gamma_r)>l_0$ then this follows from a direct calculation using the choice of $r_0$. Secondly, we use the Lipschitz condition on $u|_{S^1}$ to bound the energy $E^2(u\circ\varphi\circ\gamma_r)$ from above.
 Indeed, after possibly neglecting a zero measure set of $r$, we may assume that $$E^2(u\circ\varphi\circ (\gamma_{r})|_{\Gamma_1})  \leq r\cdot \int_{D\cap \partial B(z,r)}\mathcal{I}_+^2(\apmd u_w)\,d\hm^1(w) = r\cdot f'(r),$$ see the proof of Proposition~\ref{prop:int-hoelder-energy-filling}. Moreover, the fact that $u|_{S^1}$ is $\nu$-Lipschitz for some $\nu\geq 0$ implies  $$E^2(u\circ\varphi\circ(\gamma_r)|_{S^1\setminus\Gamma_1})\leq 2\pi \eta^2\nu^2 r^2$$ and hence 
\begin{equation}\label{eq:energy-curve-both-sides}
E^2(u\circ\varphi\circ\gamma_r)\leq r\cdot f'(r) + A'r^2
\end{equation}
 with $A'=2\pi \eta^2\nu^2$. Combining \eqref{eq:fr-qh-ef} with \eqref{eq:energy-curve-both-sides} yields \eqref{eq:diff-ineq-energy} with $\lambda=\eta^4$ and $A=\eta^4 C'M A'$. This completes the proof.
\end{proof}

\def\cprime{$'$} \def\cprime{$'$}


\begin{thebibliography}{10}

\bibitem{BBI01}
Dmitri Burago, Yuri Burago, and Sergei Ivanov.
\newblock {\em A course in metric geometry}, volume~33 of {\em Graduate Studies
  in Mathematics}.
\newblock American Mathematical Society, Providence, RI, 2001.

\bibitem{CL01}
Luca Capogna and Fang-Hua Lin.
\newblock Legendrian energy minimizers. {I}. {H}eisenberg group target.
\newblock {\em Calc. Var. Partial Differential Equations}, 12(2):145--171,
  2001.

\bibitem{Che99}
Jeff Cheeger.
\newblock Differentiability of {L}ipschitz functions on metric measure spaces.
\newblock {\em Geom. Funct. Anal.}, 9(3):428--517, 1999.

\bibitem{Chi07}
David Chiron.
\newblock On the definitions of {S}obolev and {BV} spaces into singular spaces
  and the trace problem.
\newblock {\em Commun. Contemp. Math.}, 9(4):473--513, 2007.

\bibitem{MD10}
Georgios Daskalopoulos and Chikako Mese.
\newblock Harmonic maps between singular spaces {I}.
\newblock {\em Comm. Anal. Geom.}, 18(2):257--337, 2010.

\bibitem{EG92}
Lawrence~C. Evans and Ronald~F. Gariepy.
\newblock {\em Measure theory and fine properties of functions}.
\newblock Studies in Advanced Mathematics. CRC Press, Boca Raton, FL, 1992.

\bibitem{Fug05}
Bent Fuglede.
\newblock The {D}irichlet problem for harmonic maps from {R}iemannian polyhedra
  to spaces of upper bounded curvature.
\newblock {\em Trans. Amer. Math. Soc.}, 357(2):757--792, 2005.

\bibitem{Fug08}
Bent Fuglede.
\newblock Harmonic maps from {R}iemannian polyhedra to geodesic spaces with
  curvature bounded from above.
\newblock {\em Calc. Var. Partial Differential Equations}, 31(1):99--136, 2008.

\bibitem{GT01}
David Gilbarg and Neil~S. Trudinger.
\newblock {\em Elliptic partial differential equations of second order}.
\newblock Classics in Mathematics. Springer-Verlag, Berlin, 2001.
\newblock Reprint of the 1998 edition.

\bibitem{GS92}
Mikhail Gromov and Richard Schoen.
\newblock Harmonic maps into singular spaces and {$p$}-adic superrigidity for
  lattices in groups of rank one.
\newblock {\em Inst. Hautes \'Etudes Sci. Publ. Math.}, (76):165--246, 1992.

\bibitem{Haj96}
Piotr Haj{\l}asz.
\newblock Sobolev spaces on an arbitrary metric space.
\newblock {\em Potential Anal.}, 5(4):403--415, 1996.

\bibitem{HKST15}
Juha Heinonen, Pekka Koskela, Nageswari Shanmugalingam, and Jeremy Tyson.
\newblock {\em Sobolev spaces on metric measure spaces}, volume~27 of {\em New
  Mathematical Monographs}.
\newblock Cambridge University Press, Cambridge, 2015.

\bibitem{HKST01}
Juha Heinonen, Pekka Koskela, Nageswari Shanmugalingam, and Jeremy~T. Tyson.
\newblock Sobolev classes of {B}anach space-valued functions and quasiconformal
  mappings.
\newblock {\em J. Anal. Math.}, 85:87--139, 2001.

\bibitem{Jos94}
J{\"u}rgen Jost.
\newblock Equilibrium maps between metric spaces.
\newblock {\em Calc. Var. Partial Differential Equations}, 2(2):173--204, 1994.

\bibitem{Jos97}
J{\"u}rgen Jost.
\newblock Generalized {D}irichlet forms and harmonic maps.
\newblock {\em Calc. Var. Partial Differential Equations}, 5(1):1--19, 1997.

\bibitem{Kar07}
M.~B. Karmanova.
\newblock Area and co-area formulas for mappings of the {S}obolev classes with
  values in a metric space.
\newblock {\em Sibirsk. Mat. Zh.}, 48(4):778--788, 2007.

\bibitem{Kin94}
Juha Kinnunen.
\newblock Higher integrability with weights.
\newblock {\em Ann. Acad. Sci. Fenn. Ser. A I Math.}, 19(2):355--366, 1994.

\bibitem{KS01}
Juha Kinnunen and Nageswari Shanmugalingam.
\newblock Regularity of quasi-minimizers on metric spaces.
\newblock {\em Manuscripta Math.}, 105(3):401--423, 2001.

\bibitem{Kir94}
Bernd Kirchheim.
\newblock Rectifiable metric spaces: local structure and regularity of the
  {H}ausdorff measure.
\newblock {\em Proc. Amer. Math. Soc.}, 121(1):113--123, 1994.

\bibitem{KS93}
Nicholas~J. Korevaar and Richard~M. Schoen.
\newblock Sobolev spaces and harmonic maps for metric space targets.
\newblock {\em Comm. Anal. Geom.}, 1(3-4):561--659, 1993.

\bibitem{KRS03}
Pekka Koskela, Kai Rajala, and Nageswari Shanmugalingam.
\newblock Lipschitz continuity of {C}heeger-harmonic functions in metric
  measure spaces.
\newblock {\em J. Funct. Anal.}, 202(1):147--173, 2003.

\bibitem{Lin97}
Fang~Hua Lin.
\newblock Analysis on singular spaces.
\newblock In {\em Collection of papers on geometry, analysis and mathematical
  physics}, pages 114--126. World Sci. Publ., River Edge, NJ, 1997.

\bibitem{LW15-asymptotic}
Alexander Lytchak and Stefan Wenger.
\newblock On metric spaces with a quadratic isoperimetric inequality for
  curves.
\newblock {\em in preparation}.

\bibitem{LW15-Plateau}
Alexander Lytchak and Stefan Wenger.
\newblock Area minimizing discs in metric spaces.
\newblock {\em preprint arXiv:1502.06571v2}, 2015.

\bibitem{LW-energy-area}
Alexander Lytchak and Stefan Wenger.
\newblock Energy and area minimizers in metric spaces.
\newblock {\em preprint arXiv:1507.02670v2}, 2015.

\bibitem{MZ10}
Chikako Mese and Patrick~R. Zulkowski.
\newblock The {P}lateau problem in {A}lexandrov spaces.
\newblock {\em J. Differential Geom.}, 85(2):315--356, 2010.

\bibitem{Mor48}
Charles~B. Morrey, Jr.
\newblock The problem of {P}lateau on a {R}iemannian manifold.
\newblock {\em Ann. of Math. (2)}, 49:807--851, 1948.

\bibitem{Res97}
Yu.~G. Reshetnyak.
\newblock Sobolev classes of functions with values in a metric space.
\newblock {\em Sibirsk. Mat. Zh.}, 38(3):657--675, iii--iv, 1997.

\bibitem{Res04}
Yu.~G. Reshetnyak.
\newblock Sobolev classes of functions with values in a metric space. {II}.
\newblock {\em Sibirsk. Mat. Zh.}, 45(4):855--870, 2004.

\bibitem{Res06}
Yu.~G. Reshetnyak.
\newblock On the theory of {S}obolev classes of functions with values in a
  metric space.
\newblock {\em Sibirsk. Mat. Zh.}, 47(1):146--168, 2006.

\bibitem{SU82}
Richard Schoen and Karen Uhlenbeck.
\newblock A regularity theory for harmonic maps.
\newblock {\em J. Differential Geom.}, 17(2):307--335, 1982.

\bibitem{Ser94}
Tomasz Serbinowski.
\newblock Boundary regularity of harmonic maps to nonpositively curved metric
  spaces.
\newblock {\em Comm. Anal. Geom.}, 2(1):139--153, 1994.

\bibitem{Str80}
Edward~W. Stredulinsky.
\newblock Higher integrability from reverse {H}\"older inequalities.
\newblock {\em Indiana Univ. Math. J.}, 29(3):407--413, 1980.

\bibitem{Tuk80}
Pekka Tukia.
\newblock The planar {S}ch\"onflies theorem for {L}ipschitz maps.
\newblock {\em Ann. Acad. Sci. Fenn. Ser. A I Math.}, 5(1):49--72, 1980.

\bibitem{Wen07}
Stefan Wenger.
\newblock Flat convergence for integral currents in metric spaces.
\newblock {\em Calc. Var. Partial Differential Equations}, 28(2):139--160,
  2007.

\bibitem{Zap14}
Aleksandra Zapadinskaya.
\newblock H\"older continuous {S}obolev mappings and the {L}usin {N} property.
\newblock {\em Illinois J. Math.}, 58(2):585--591, 2014.

\bibitem{ZZ14}
Hui-Chun Zhang and Xi-Ping Zhu.
\newblock Lipschitz continuity of harmonic maps between alexandrov spaces.
\newblock {\em preprint arXiv:1311.1331}, 2013.

\end{thebibliography}
\end{document}